\documentclass[11pt,reqno]{amsart}
\usepackage{amsmath,amsfonts, amssymb, amsthm}
\usepackage{mathrsfs}
\usepackage{graphicx, color}
\usepackage[numbers]{natbib} 
\usepackage[colorlinks,citecolor=green]{hyperref}
\usepackage{tikz}

\def\qed{\hfill \rule{4pt}{7pt}}
\def\qed{\hfill \rule{4pt}{7pt}}
\begin{document}
	
	\title[The combinatorics of some two-color partition identities  ]
	{The combinatorics of some two-color partition identities}
	
	
	\author[Y.-C. Shen]{Yong-Chao Shen}
	\address{School of Mathematics and KL-AAGDM, Tianjin
		University, Tianjin 300350, People's Republic of China}
	\email{{\tt yc\_s@tju.edu.cn}}

	
	\keywords{two color partition, combinatorial  proofs 
		\newline \indent 2020 {\it Mathematics Subject Classification}. Primary 05A19, 11B65; Secondary 33B15.
		\newline \indent Supported by National Key Research and Development Program of China
		(2023YFA1009401).}
	
	\begin{abstract} 
		Recently, Andrews and EI Bachraoui obtained several identities on two-colored partitions. While solving open problems they posed, Chen and Zhou derived a number of identities using analytic methods and asked for combinatorial proofs. In this note, we provide the requested combinatorial proofs. Additionally, we derive several new identities and provide combinatorial proofs for them.
		
	\end{abstract}
	\maketitle

	\section{Introduction}
	
	In recent years, the study of two-colored partitions has attracted the attention of numerous scholars. Andrews and El Bachraoui \citep{1,2,3, 4}  have given numerous identities and inequalities on two-colored partitions analytically. Chen and Liu \cite{8}, Chen and Zou \cite{9}, Bugleev \cite{12}, Fu \cite{10}, Gu and Yu \cite{11} provide combinatorial proofs for several identities on two-colored partitions.
	
	An overpartition of $n$ is a partition of $n$ where the first occurrence of each part may be overlined and  denote the cardinality of it by $\overline{p}(n)$. For example, $(\overline{5}, 5,3,2,\overline{1},1)$ is an overpartition of $17$. We denote the number of overpartitions of $n$ into odd parts by $\overline{p}_o(n)$. 
	
	Recently, Andrews and El Bachraoui \cite{1} gave some identities on two color partition by using analytical methods and  Chen and Zhou \cite{9} analogously derived several identities using analytical methods. They introduced $F(n)$ as the number of two-colour partitions of $n$ in which  the even parts may occur only in the blue colour. They denoted by $F_0(n)$ and $F_1(n)$ the number of partitions counted by $F(n)$ in which the number of even parts is respectively even or odd, and $F_2(n)$ and $F_3(n)$ as the
	number of partitions counted by $F(n)$ in which the number of parts is respectively even or odd. 
	
	They derived the following identities.

	


\textbf{ Theorem 1.1.} (\cite{9} Theorem 1.3.) For any nonnegative integer $n$,

\medskip

(a) $F(n) = \overline{p}(n)$,

(b) $F_0(n) = \frac{1}{2}\bigl(\overline{p}(n) + \overline{p}_o(n)\bigr)$,

(c) $F_1(n) = \frac{1}{2}\bigl(\overline{p}(n) - \overline{p}_o(n)\bigr)$,

(d) $F_2(n) = \frac{1}{2}\bigl(\overline{p}(n) + (-1)^n \overline{p}_o(n)\bigr)$,

(e) $F_3(n) = \frac{1}{2}\bigl(\overline{p}(n) - (-1)^n \overline{p}_o(n)\bigr)$.

Furthermore, we have also arrived at similar results.

\textbf{ Definition 1.2.}  Let $n$ be a nonnegative integer and let $G(n) $be the number of two color partitions of $n$ in which  the odd parts may occur only in the blue color. Let $G_0(n)$(resp. $G_1(n)$) be the number of partitions counted by $G(n)$ in which the number of blue even parts is even (resp. odd). Furthermore, let $G_2(n)$ (resp. $G_3(n)$) be the number of partitions counted by $G(n)$ in which the number of blue parts is even (resp. odd). Let $G_4(n)$ (resp. $G_5(n)$) be the number of partitions counted by $G(n)$ in which the number of  parts is even (resp. odd).

\textbf{ Definition 1.3.}  Let $n$ be a nonnegative integer and let $H(n) $be the number of  partitions of $n$ in which  the odd parts are distinct.

\textbf{ Definition 1.4.}
Let $n$ be a nonnegative integer and let $K(n) $be the number of  partitions of $n$ in which 

$\bullet$ each part is distinct;

$\bullet$ even parts can only be congruent to 2 modulo 4;

$\bullet$ the count is weighted by $(-1)^j$ where $j$ is the number of the even parts.

It is easy to see that\[
\sum_{n=0}^{\infty} K(n)q^n=(q^2;q^4)_\infty (-q;q^2)_\infty=\frac{(-q;q^2)_\infty}{(-q^2;q^2)_\infty}.
\]
\textbf{ Theorem 1.5.} For any nonnegative integer $n$,

$(a')$ $G_0(n) = \frac{1}{2}\bigl(G(n) + H(n)\bigr)$,

$(b')$ $G_1(n) = \frac{1}{2}\bigl(G(n) - H(n)\bigr)$,

$(c')$ $G_2(n) = \frac{1}{2}\bigl(G(n) +(-1)^n H(n)\bigr)$,

$(d')$ $G_3(n) = \frac{1}{2}\bigl(G(n) -(-1)^n H(n)\bigr)$,

$(e')$ $G_4(n) = \frac{1}{2}\bigl(G(n) +(-1)^n K(n)\bigr)$,

$(f')$ $G_5(n) = \frac{1}{2}\bigl(G(n) -(-1)^n K(n)\bigr)$.

\textbf{ Definition 1.6.} Let $n$ be a nonnegative integer and let $L(n) $be the number of two color partitions of $n$ in which all parts are distinct and the odd parts may occur only in the blue color. Let $L_0(n)$(resp. $L_1(n)$) be the number of partitions counted by $L(n)$ in which the number of blue even parts is even (resp. odd). Furthermore, let $L_2(n)$ (resp. $L_3(n)$) be the number of partitions counted by $L(n)$ in which the number of blue parts is even (resp. odd).

\textbf{ Theorem 1.7.} For any nonnegative integer $n$, there holds

(f) \[
L_0(n) - L_1(n) = 
\begin{cases} 
	1 & \text{if }  n  \text{ is a }   \text{triangular number,}\\
	0 & \text{otherwise}.
\end{cases}
\]

(g)\[
L_2(n) - L_3(n) = 
\begin{cases} 
	(-1)^n & \text{if }  n \text{ is a  }   \text{triangular number,}\\
	0 & \text{otherwise}.
\end{cases}
\]

For (a), Chen and Zhou have gave a combinatorial proof of it. In Section 2, we provide a combinatorial proof of  (b),(c),(d),(e),$(a') $, $(b') $,$(c')$ and $(d')$and provide analytic proof of $(e')$ and $(f')$. In Section 3, we provide the proof of $(f) $ and $(g)$.

\section{ The  proof of $(b)$,$(c)$,$(d)$,$(e)$,$(a') $, $(b')$ $(c')$ , $(d')$, $(e')$and $(f')$}

	We denote $\pi_{\overline{p}_o}(n)$ the set of overpartitions of $n$ into odd parts and $\pi _F(n)$ the set of two-colour partitions of $n$ in which the even parts may occur only in the blue colour. We also denote $\pi _G(n)$ the set of two-colour partitions of $n$ in which the odd parts may occur only in the blue colour and $\pi _H(n)$ the set of partitions of $n$ in which the odd parts are distinct.
	
	$Proof$ $ of$ $ (b)$ and $(c).$
	It suffices to prove that  

\[	F_0(n)-F_1(n)=\overline{p}_o(n).\]
	
	Let $\pi_ Q(n)$  be a subset of $\pi_ F(n)$, where the blue parts are distinct odd numbers. $Q(n)$ counts the number of $\pi_ Q(n)$. We can define an involution $\Phi$ which changes the parity of the number of even parts on $\pi_ F(n)$  $\setminus$ $\pi_ Q(n)$ :
	
	For any $\lambda$ $\in$ $\pi_ F(n)$  $\setminus$ $\pi_ Q(n)$ .
	
	(1) If there is no even number in $\lambda$, then  there must exists odd numbers of equal size in the blue color in the $\lambda$. Let 
	$a$	denotes the largest such odd number.  $\mu$ is the partition obtained from $\lambda$	by merging the two $a$ into $2a$. Then $\mu$ $\in$ $\pi_ F(n)$  $\setminus$ $\pi_ Q(n)$ .

	(2) If there has even numbers in $\lambda$, we denote the largest even number by $c$. We have two subcases.

	(2-1) If there are odd numbers of equal size in the blue color, denote the largest such number by $d$.
	
	(2-1-i) If $2d > c$,  then merging the two $d$ into $2d$ to obtain $\mu$. Then $\mu$ $\in$ $\pi_ F(n)$  $\setminus$ $\pi_ Q(n)$ .
	
(2-1-ii) When the largest repeated even part is not $c$ in $\lambda$.

	If $2d \leq c$ and $c=4k+2$, dividing $c$ into two equal parts to obtain $\mu$. Then $\mu$ $\in$ $\pi_ F(n)$  $\setminus$ $\pi_ Q(n)$ .
	
	When $2d \leq c$ and $c=4k$, we have two subcases.
	
	If the other even parts in $\lambda$ do not exceed $\frac{c}{2}$ or the even parts exceeding $\frac{c}{2}$ are pairwise distinct, dividing $c$ into two equal parts to obtain $\mu$. Then $\mu$ $\in$ $\pi_ F(n)$  $\setminus$ $\pi_ Q(n)$ .
	

	If there exist repeated even parts exceeding $\frac{c}{2}$, let the maximal such even number be denoted as e. Then merging the two $e$ into $2e$ to obtain $\mu$. Then $\mu$ $\in$ $\pi_ F(n)$  $\setminus$ $\pi_ Q(n)$ .
	
	
(2-1-iii)	If $2d \leq c$ and the largest repeated even part is  $c$ in $\lambda$,  then merging the two $c$ into $2c$ to obtain $\mu$. Then $\mu$ $\in$ $\pi_ F(n)$  $\setminus$ $\pi_ Q(n)$ .
	
	(2-2) When there no odd number of equal size in the blue color.
	
	(2-2-i) If the largest repeated even part is  $c$ in $\lambda$,   then merging the two $c$ into $2c$ to obtain $\mu$. Then $\mu$ $\in$ $\pi_ F(n)$  $\setminus$ $\pi_ Q(n)$ .
	
	(2-2-ii)When the largest repeated even part is not $c$ in $\lambda$.
	
	When $c=4k+2$, 	then dividing the largest even number $c$ into two equal parts to obtain $\mu$. Then $\mu$ $\in$ $\pi_ F(n)$  $\setminus$ $\pi_ Q(n)$.
	
	When $c=4k$, we have two subcases.
	 
	 If the other even parts in $\lambda$ do not exceed $\frac{c}{2}$ or the even parts exceeding $\frac{c}{2}$ are pairwise distinct, dividing $c$ into two equal parts to obtain $\mu$. Then $\mu$ $\in$ $\pi_ F(n)$  $\setminus$ $\pi_ Q(n)$ .
	
	
		If there exist repeated even parts exceeding $\frac{c}{2}$, let the maximal such even number be denoted as $f$. Then merging the two $f$ into $2f$ to obtain $\mu$. Then $\mu$ $\in$ $\pi_ F(n)$  $\setminus$ $\pi_ Q(n)$ .
		

\textbf{Example 2.1.}   We give some examples as follows. 
		
	\begin{tabular}{c|c} 
		$8_b+1_b$ & $4_b+4_b+1_b$ \\
	$5_b+2_b+1_b+1_g$	 & $5_b+1_b+1_b+1_b+1_g$  \\
	$8_b+8_b+3_b+3_b+3_b+1_b+1_b+1_g$ & $16_b+3_b+3_b+3_b+1_b+1_b+1_g$ \\
	$7_g+6_b+4_b+4_b+3_b+1_b+1_b$ & $7_g+4_b+4_b+3_b+3_b+3_b+1_b+1_b$  \\
	\end{tabular}
	
	\vspace{1cm} 
	
Futhermore,	we construct a map $\phi$: $\pi_ Q(n)$  $\rightarrow$ $\pi_{\overline{p}_o}(n)$ as follows.  
	
	For any $\lambda$ in $\pi_ Q(n)$, map the green parts in $\lambda$ to the nonoverlined parts and map the blue parts in $\lambda$ to the overlined parts to obtain the $\beta$=$\phi$($\lambda$). Then $\beta$ $\in$ $\pi_{\overline{p}_o}(n)$.
	
 Obviously, $\phi$ is a bijection. This implies that $Q(n) =\overline{p}_o(n).$ 
 
 Therefore
 \[	F_0(n)-F_1(n)=\overline{p}_o(n).\]\qed
 
$ Proof$ $of $ $(d)$ \text{and} $(e).$
 It suffices to prove that
 
 \[	F_2(n)-F_3(n)=(-1)^n \overline{p}_o(n).\] 
 
 If $n$ be even then it has two subcases.
 
 (1) The number of even parts is even and the number of odd parts is even, then the number of the parts is even;
 
 (2)The number of even parts is odd and the number of odd parts is even, then  the number of the parts is odd;
 
 therefore
 \[	F_2(n)-F_3(n)=F_0(n)-F_1(n)=\overline{p}_o(n).\]
 
  If $n$ be odd then it has two subcases.
 
 (1) The number of even parts is even and the number of odd parts is odd, then the number of the parts is odd;
 
 (2)The number of even parts is odd and the number of odd parts is odd, then  the number of the parts is even;
 
 therefore
 \[	F_2(n)-F_3(n)=F_1(n)-F_0(n)=-\overline{p}_o(n).\]
 
 So \[	F_2(n)-F_3(n)=(-1)^n \overline{p}_o(n).\]\qed
 
 $ Analytic$  $ proof
 $ $ of$ $ (a')$ , $(b')$, $(c')$ and $(d')$.
 We can easily to see 
 
 \[
 \sum_{n=0}^{\infty} G(n) q^n=\frac{1}{(q^2;q^2)_{\infty}(q,q)_{\infty}},
 \]
 and
 \[
 \sum_{n=0}^{\infty} H(n)q^n=\frac{(-q;q^2)_\infty}{(q^2;q^2)_\infty}.
 \]
 
 Then
 \[
 \sum_{n=0}^{\infty} (G_0(n)+G_1(n)) q^n=\sum_{n=0}^{\infty} (G_2(n)+G_3(n)) q^n=\sum_{n=0}^{\infty} G(n)q^n,
 \]
 and 
 \[
 \sum_{n=0}^{\infty} (G_0(n)-G_1(n)) q^n=\frac{1}{(-q^2;q^2)_\infty(q;q)_\infty}=\frac{(-q;q^2)_\infty}{(q^2,q^2)_\infty}=\sum_{n=0}^{\infty} H(n)q^n,
 \]
  and 
  \[
  \sum_{n=0}^{\infty} (G_2(n)-G_3(n)) q^n=\frac{1}{(q^2;q^2)_\infty(-q;q)_\infty}=\sum_{n=0}^{\infty}(-1)^n H(n)q^n.
  \]

 Combining these three formulae, we can complete the proof.\qed

  $Combinatorial$   $ proof$
  $ of$ $ (a')$ and $(b').$
 It suffices to prove that  
 
 \[	G_0(n)-G_1(n)=H(n).\]
 
 Let $\pi_ R(n)$  be a subset of $\pi_ G(n)$, which the blue parts are distinct odd numbers. $R(n)$ counts the number of $\pi_ R(n)$. We can use the involution  which is used to prove (b) and (e)  on $\pi_ G(n)$  $\setminus$ $\pi_ R(n)$ .
 
 Therefore, we only need to construct a bijection $\psi$ : $\pi_ R(n)$ $\rightarrow$ $\pi_ H(n)$ to complete the proof.
 
 For any $\lambda$ $\in$ $\pi_ R(n)$, clearing the all color to obtain the $\mu$, then $\mu$ $\in$ $\pi_ H(n)$.
 
  Conversely,  For any $\mu$ $\in$ $\pi_ H(n)$, marking all the odd numbers in blue and marking all the even numbers in green to obtain $\lambda$, then $\lambda$ $\in$ $\pi_ R(n)$.
  
   Obviously, $\psi$ is a bijection. This implies that $R(n) =H(n).$ 
  
  Therefore
  \[	G_0(n)-G_1(n)=H(n).\]\qed
  
  $Combinatorial$   $ proof$
  $ of$ $ (c')$ and $(d').$
  It suffices to prove that  
  
  \[	G_2(n)-G_3(n)=(-1)^nH(n).\]
  
  For any $\lambda$ $\in$ $\pi_G(n)$, the sum of green parts is even.
  
   If $n$ be even then the sum of blue parts is even. It has two subcases.
  
  (1) The number of blue even parts is even and the number of blue odd parts is even, then the number of the blue parts is even;
  
  (2)The number of blue even parts is odd and the number of blue odd parts is even, then  the number of the blue parts is odd;
  
  therefore
  \[	G_2(n)-G_3(n)=G_0(n)-G_1(n)=H(n).\]
  
  If $n$ be odd then the sum of blue parts is odd.  It has two subcases.
  
  (1) The number of blue even parts is even and the number of blue odd parts is odd, then the number of blue the parts is odd;
  
  (2)The number of blue even parts is odd and the number of blue odd parts is odd, then  the number of the blue parts is even;
  
  therefore
  \[	G_2(n)-G_3(n)=G_1(n)-G_0(n)=-H(n).\]
  
  So \[	G_2(n)-G_3(n)=(-1)^n H(n).\]\qed

   $ Analytic$  $ proof
  $ $ of$ $ (e')$ and $(f').$ We can get
  \[
  \sum_{n=0}^{\infty} (G_4(n)-G_5(n)) q^n=\frac{1}{(-q^2;q^2)_\infty(-q;q)_\infty}=\frac{(q;q^2)_\infty}{(-q^2,q^2)_\infty}=\sum_{n=0}^{\infty}(-1)^n K(n)q^n,
  \]
  
  and 
  \[
  \sum_{n=0}^{\infty} (G_4(n)+G_5(n)) q^n=\sum_{n=0}^{\infty} G(n)q^n.
  \]
  
  Combining these two formulae, we can complete the proof.\qed
  
  It is natural to ask if there are combinatorial proofs for parts $(e')$ and $(f')$.
 
 \section{ The  proof of $(f)$ and $(g)$}
 
 $ Analytic$  $ proof
 $ $ of$ $ (f)$ and $(g).$
 We can easily to see 
 \[
 \sum_{n=0}^{\infty} L(n)q^n=(-q^2,q^2)_\infty(-q,q)_\infty.
 \]
 
 Then\[
 \sum_{n=0}^{\infty} (L_0(n)-L_1(n))q^n=(q^2,q^2)_\infty(-q,q)_\infty=\frac{(q^2,q^2)_\infty}{(q,q^2)_\infty}=\sum_{n=0}^{\infty}q^{\frac{n(n+1)}{2}},
 \]
 and
 \[
 \sum_{n=0}^{\infty} (L_2(n)-L_3(n))q^n=(-q^2,q^2)_\infty(q,q)_\infty=\sum_{n=0}^{\infty}(-1)^nq^{\frac{n(n+1)}{2}}. 
 \] 
 Combining these two formulae, we can complete the proof.\qed

 Before providing a combinatorial proof for $(f)$ and $(g)$, we first introduce the $4$-modular diagram of $\lambda$ which all parts are distinct and the even part must be a multiple of $4$. For $\lambda$, we can decompose it into three parts, $\lambda_e$ which parts are even, $\lambda_{C1}$ which parts are congruent to $1$ modulo $4$, and $\lambda_{C3}$ which parts are congruent to $3$ modulo $4$. In the $4$-modular diagram, $\lambda_{C3}$ parts are represented in rows above the main diagonal and $\lambda_{C1}$ parts are represented in columns below the main diagonal ( $\lambda_e$ parts are not represented in the diagram). For $4k+1$, its representation would have $k$ squares and $1$ triangle and write $1$ inside the triangle. For $4k+3$, its representation would have $k$ squares and $1$ triangle and write $3$ inside the triangle. 
 
 We give two examples as follows.
 
 \textbf{Example 3.1.}  $\lambda_e$=$(8,4)$, $\lambda_{C1}$=$(5,1)$, $\lambda_{C3}$=$(3)$.
 
\begin{center}
\begin{tikzpicture}
	\draw (0,0) -- (0,2);
	\draw (0,0) -- (2,0);
	\draw (2,0) -- (2,1);
	\draw (1,0) -- (1,1);
	\draw (0,1) -- (2,1);
	\draw (0,2) -- (2,0);
	
	\node at (0.2, 1.5) {1};
	\node at (1.2, 0.5) {1};
	\node at (1.8, 0.6) {3};
\end{tikzpicture}
\end{center}

 \textbf{Example 3.2.}  $\lambda_e$=$(12,8,4)$, $\lambda_{C1}$=$(17,13,9,5,1)$, $\lambda_{C3}$=$(19,15,11,7,3)$.
 
 
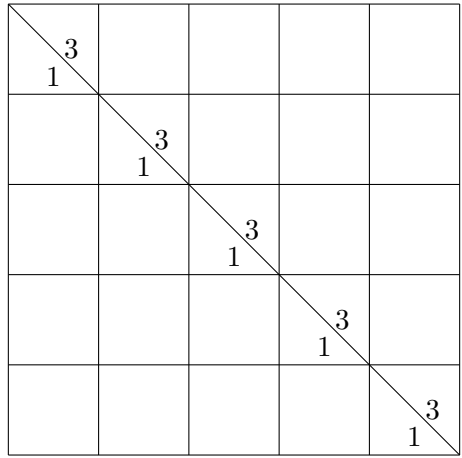
\begin{center}
 \begin{tikzpicture}[scale=1.2]
 	\foreach \x in {0,...,5} {
 		\draw (\x,0) -- (\x,5);
 		\draw (0,\x) -- (5,\x);
 	}
 	\draw (0,5) -- (5,0);
 	\node at (0.5, 4.2) {1};
 	\node at (1.5, 3.2) {1};
 	\node at (2.5, 2.2) {1};
 	\node at (3.5, 1.2) {1};
 	\node at (4.5, 0.2) {1};
 	\node at (0.7, 4.5) {3};
 	\node at (1.7, 3.5) {3};
 	\node at (2.7, 2.5) {3};
 	\node at (3.7, 1.5) {3};
 	\node at (4.7, 0.5) {3};
 \end{tikzpicture}
 	\end{center}
 	
 $Combinatorial$ $proof$  $for$ $(f)$ $and$ $(g)$. For $(f)$, we let $\pi _L(n)$ the set of two-colour partitions of $n$ in which all parts are distinct and the odd parts may occur only in the blue colour. We also denote $\pi _M(n)$ the subset of $\pi _L(n)$ in which there exists at least an even number in only one color.
 Also, we let $\pi _N(n)$ the set of  partitions of $n$ in which all parts are distinct and the even part must be a multiple of $4$.

 Firstly, we can construct an involution $\theta$ that changes the sign on $\pi _M(n)$: change the color of the largest even number that is only in one color. Obviously, $\theta$ changes the parity of the number of even numbers in the blue part. For example, 
 $(5_b,4_b,3_b,2_g,2_b)$ $\rightarrow$ $(5_b,4_g,3_b,2_g,2_b)$.
 
 Next, we only to consider on  $\pi _L(n)$ $\setminus$  $\pi _M(n)$. By definition, we know that the same even number appears once in each of the two colors in $\lambda$ which is in  $\pi _L(n)$ $\setminus$  $\pi _M(n)$. 
 
 For any $\lambda$ $\in$ $\pi _L(n)$ $\setminus$  $\pi _M(n)$, merge even numbers of the same size and clear the color of all parts to obtain the $\mu$. Then $\mu$ $\in$  $\pi _N(n)$. Conversely, for any $\mu$ $\in$ $\pi _N(n)$, bisect each even number, label the odd numbers blue, and for even numbers of the same magnitude, label one blue and the other green to obtain $\lambda$. Then $\lambda$ $\in$ $\pi _L(n)$ $\setminus$  $\pi _M(n)$. This indicates that there exists a bijection between $\pi _L(n)$ $\setminus$  $\pi _M(n)$ and $\pi _N(n)$. Also, the number of blue even parts of $\lambda$ equals to the number of even parts of $\mu$.
 
Consider the 4-modular diagram of $\lambda$ which is in $\pi _N(n)$.  Let us define a transformation from $\pi _N(n)$ onto itself:
 
 1.If there are two triangles adjoined, forming a square with diagonal line, then delete the diagonal line.
 
 2. The next step of the transformation depends on the cardinality of $\lambda_{C_1}$ and $\lambda_{C3}$. 
 
 2.1 If $|\lambda_{C3}|$ $>$ $|\lambda_{C_1}|$, then let $l$ denote the length of the longest column containing only squares. If product $4l$ is strictly greater than the maximum part in $\lambda_e$ (let us agree that max $\lambda_e$ = $0$, if $\lambda_e$ = $\emptyset$), then remove that longest column and include $4l$ as a part in $\lambda_e$. Otherwise, place the max $\lambda_e$ part into the diagram as a column to the left of the longest column. If $l = 0$ and $\lambda_e$ $\neq$ $\emptyset$ then the max $\lambda_e$ part should be added to the far right of the diagram as a column.
 
 2.2 If $|\lambda_{C3}|$ $\leq$ $|\lambda_{C_1}|$ then let $l$ denote the longest row containing only squares. If product $4l$ is strictly greater than max $\lambda_e$, then remove that row and include $4l$ as a part in $\lambda_e$. Otherwise, place the max $\lambda_e$ part as a row into the diagram right above the longest row of squares. If $l = 0$ and $\lambda_e$ $\neq$ $\emptyset$ then max $\lambda_e$ part should be added at the bottom of the diagram as a row.
 
 3. After that draw back the diagonal line and read off the new partition $\lambda$=($\lambda'_e$,$\lambda'_{C3}$,$\lambda'_{C1}$).
 
 There are two partitions for which the transformation is not applicable: $\lambda = ((4k-3), . . . ,9,5,1)$ and $\lambda = ((4k-1), . . . ,11,7,3)$. Here the mapping does not work, because such partitions do not have any rows or columns containing only squares and $\lambda_e$ = $\emptyset$. However, for all other partitions the transformation establishes an involution.
 
  For $\lambda = ((4k-3), . . . ,9,5,1)$, the weight of $\lambda$ is $ 2k^2-k$, We know that triangular number $T_k=\frac{k(k+1)}{2}$, let $k=2n-1$, then $T_{2n-1}=2n^2-n$, so the weight of $\lambda$ is $T_{2n-1}$.
  
  For $\lambda = ((4k-1), . . . ,11,7,3)$, the weight of $\lambda$ is $ 2k^2+k$, We know that triangular number $T_k=\frac{k(k+1)}{2}$, let $k=2n$, then $T_{2n}=2n^2+n$, so the weight of $\lambda$ is $T_{2n}$.
 
  We note that when $k=0, 2k^2+k=2k^2-k$. Therefore for $k>0$, there is only one partition of $T_k$ for which the transformation is not applicable.
 
 Therefore, for $\pi _L(n)$ , $\lambda = ((4k-3)_b, . . . ,9_b,5_b,1_b)$ and $\lambda = ((4k-1)_b, . . . ,11_b,7_b,3_b)$ can not be transformed.
 
 Hence the $(f)$ is proved.
 
 It deserves to be mentioned, that the same idea of an almost-involution happens to be in the proof of Euler’s pentagonal theorem \cite{5}. Moreover, the generalization of 4-modular diagram allows to prove combinatorically the Jacobi triple product identity \cite{15}  and the generalization of 2-modular diagram is  used in \cite{12}.

 For $(g)$, $\theta$ also changes the parity of the number of blue parts in $\pi_M(n). $ We only need to pay attention to $\pi _L(n)$ $\setminus$  $\pi _M(n)$. By definition, for any $\lambda$ $\in$ $\pi _L(n)$ $\setminus$  $\pi _M(n)$, the sum of green parts of $\lambda$ is even.
 
 If $n$ be even  triangular number then the sum of blue parts of $\lambda$ is even.  It has two subcases.
 
 (1) The number of blue even parts is even and the number of blue odd parts is even, then the number of the blue parts is even;
 
 (2)The number of blue even parts is odd and the number of blue odd parts is even, then  the number of the blue parts is odd;
 
 therefore
 \[	L_2(n)-L_3(n)=L_0(n)-L_1(n)=1.\]
 
 If $n$ be odd  triangular number then the sum of blue parts of $\lambda$ is odd. It has two subcases.
 
 (1) The number of blue even parts is even and the number of blue odd parts is odd, then the number of the blue parts is odd;
 
 (2)The number of blue even parts is odd and the number of blue odd parts is odd, then  the number of the blue parts is even;
 
 therefore
 \[	L_2(n)-L_3(n)=L_1(n)-L_0(n)=-1.\]
 
 So \[	L_2(n)-L_3(n)=(-1)^n .\]\qed

 \textbf{Example 3.3.} Let us give an example.
 
 $\gamma$=$(6_g,6_b,5_b,4_g,4_b,3_b,2_g,2_b,1_b)$ $\in$ $\pi _L(n)$, then $\lambda$=$(12,8,5,4,3,1)$ $\in$ $\pi _N(n)$.
 
 $\lambda_e$=$(12,8,4)$, $\lambda_{C1}$=$(5,1)$, $\lambda_{C3}$=$(3)$.

 \begin{tikzpicture}
 	\draw (0,0) -- (0,2);
 	\draw (0,0) -- (2,0);
 	\draw (2,0) -- (2,1);
 	\draw (1,0) -- (1,1);
 	\draw (0,1) -- (2,1);
 	\draw (0,2) -- (2,0);
 	
 	\node at (0.2, 1.5) {1};
 	\node at (1.2, 0.5) {1};
 	\node at (1.8, 0.6) {3};
 \end{tikzpicture} $\rightarrow$
 \begin{tikzpicture}
 	\draw (0,0) -- (0,2);
 	\draw (0,0) -- (2,0);
 	\draw (2,0) -- (2,1);
 	\draw (1,0) -- (1,1);
 	\draw (0,1) -- (2,1);
 	\draw (0,2) -- (1,1);
 	
 \end{tikzpicture} $\rightarrow$
 \begin{tikzpicture}
 	\draw (0,0) -- (0,3);
 	\draw (0,0) -- (2,0);
 	\draw (2,0) -- (2,1);
 	\draw (1,0) -- (1,2);
 	\draw (0,1) -- (3,1);
 	\draw (0,2) -- (3,2);
 	\draw (2,0) -- (2,2);
 	\draw (3,1) -- (3,2);
 	\draw (0,3) -- (1,2);
 \end{tikzpicture} $\rightarrow$
 \begin{tikzpicture}
 	\draw (0,0) -- (0,3);
 	\draw (0,0) -- (2,0);
 	\draw (2,0) -- (2,1);
 	\draw (1,0) -- (1,2);
 	\draw (0,1) -- (3,1);
 	\draw (0,2) -- (3,2);
 	\draw (2,0) -- (2,2);
 	\draw (3,1) -- (3,2);
 	\draw (0,3) -- (2,1);
 	\node at (0.2, 2.2) {1};
 	\node at (1.2, 1.2) {1};
 		\node at (1.8, 1.4) {3};
 \end{tikzpicture}
 
 Then $\lambda'_e$=$(8,4)$, $\lambda'_{C1}$=$(9,5)$, $\lambda'_{C3}$=$(7)$.
 $\gamma'$=$(9_b,7_b,5_b,4_g,4_b,2_g,2_b)$$\in$ $\pi _L(n)$.

\end{document}